\input amstex
\documentstyle{amsppt}
%
\catcode`@=11
\redefine\output@{%
  \def\break{\penalty-\@M}\let\par\endgraf
  \ifodd\pageno\global\hoffset=105pt\else\global\hoffset=8pt\fi  
  \shipout\vbox{%
    \ifplain@
      \let\makeheadline\relax \let\makefootline\relax
    \else
      \iffirstpage@ \global\firstpage@false
        \let\rightheadline\frheadline
        \let\leftheadline\flheadline
      \else
        \ifrunheads@ 
        \else \let\makeheadline\relax
        \fi
      \fi
    \fi
    \makeheadline \pagebody \makefootline}%
  \advancepageno \ifnum\outputpenalty>-\@MM\else\dosupereject\fi
}
\catcode`\@=\active
\nopagenumbers
\def\negskp{\hskip -2pt}
\def\sign{\operatorname{sign}}
\def\tr{\operatorname{tr}}
\def\Rea{\operatorname{Re}}
\def\Img{\operatorname{Im}}
\def\vtrule{\vrule height 12pt depth 6pt}

\def\blue#1{#1}
\catcode`#=11\def\diez{#}\catcode`#=6
\catcode`_=11\def\podcherkivanie{_}\catcode`_=8
\def\mycite#1{\cite{\blue{#1}}\immediate\special{ps:
     ShrHPSdict begin /ShrBORDERthickness 0 def}}

\def\mytag#1{%
    \tag#1}
\def\mythetag#1{\thetag{\blue{#1}}\immediate\special{ps:
     ShrHPSdict begin /ShrBORDERthickness 0 def}}
\def\myrefno#1{\no#1}
\def\myhref#1#2{\blue{#2}\immediate\special{ps:
     ShrHPSdict begin /ShrBORDERthickness 0 def}}
\def\myEarXivlink{\myhref{http://arXiv.org}{http:/\negskp/arXiv.org}}
\def\myGeoCities{\myhref{http://www.geocities.com}{GeoCities}}
\def\mytheorem#1{\csname proclaim\endcsname{Theorem #1}}
\def\mythetheorem#1{\blue{#1}\immediate\special{ps:
     ShrHPSdict begin /ShrBORDERthickness 0 def}}
\def\mylemma#1{\csname proclaim\endcsname{Lemma #1}}

\def\mycorollary#1{\csname proclaim\endcsname{Corollary #1}}

\pagewidth{360pt}
\pageheight{606pt}
\topmatter
\title
A note on pairs of metrics\\ in a two-dimensional 
linear vector space.
\endtitle
\author
R.~A.~Sharipov
\endauthor
\address 5 Rabochaya street, 450003 Ufa, Russia\newline
\vphantom{a}\kern 12pt Cell Phone: +7(917)476 93 48
\endaddress
\email \vtop to 30pt{\hsize=280pt\noindent
\myhref{mailto:r-sharipov\@mail.ru}
{r-sharipov\@mail.ru}\newline
\myhref{mailto:R\podcherkivanie Sharipov\@ic.bashedu.ru}
{R\_\hskip 1pt Sharipov\@ic.bashedu.ru}\vss}
\endemail
\urladdr
\vtop to 20pt{\hsize=280pt\noindent
\myhref{http://www.geocities.com/r-sharipov}
{http:/\negskp/www.geocities.com/r-sharipov}\newline
\myhref{http://www.freetextbooks.boom.ru/index.html}
{http:/\negskp/www.freetextbooks.boom.ru/index.html}\vss}
\endurladdr
\abstract
    Pairs of metrics in a two-dimensional linear vector space 
are considered, one of which is a Minkowski type metric. Their 
simultaneous diagonalizability is studied and canonical 
presentations for them are suggested.
\endabstract
\endtopmatter
\loadbold
\loadeufb
\TagsOnRight
\document

\rightheadtext{A note on pairs of metrics \dots}
\head
1. Introduction. 
\endhead
    Let $V$ be a two-dimensional linear vector space equipped 
with two metrics $\bold g$ and $\check{\bold g}$. It is well 
known (see \mycite{1}) that if $\bold g$ is positive, then these
two metrics can be diagonalized simultaneously in some basis. 
Here we consider a different case where $\bold g$ is a metric
of the signature $(+,-)$. A two-dimensional space with such a
metric is often used as a two-dimensional model of the four-dimensional
Minkowski space. For this reason we call $\bold g$ a Minkowski
type metric.
\head
2. Lorentz transformations and diagonalizability.
\endhead
     Each single metric can be diagonalized. This fact means that
there is some basis $\bold e_0,\,\bold e_1$ in $V$ such that the 
metric $\bold g$ is given by the diagonal matrix $g_{ij}$ with 
the numbers $1$ and $-1$ on its diagonal:
$$
\xalignat 2
&\hskip -2em
g_{ij}=\Vmatrix\format\l\ &\ \r\\1 & 0\\ 0 & -1\endVmatrix,
&&\check g_{ij}=\Vmatrix\check g_{00} & \check g_{01}\\
\vspace{1ex}\check g_{01} & \check g_{11}\endVmatrix.
\mytag{2.1}
\endxalignat
$$
The second metric is not necessarily diagonal in this basis. In order
to diagonalize it we perform the following Lorentz transformation of
the basis $\bold e_0,\,\bold e_1$:
$$
\hskip -2em
\aligned
&\tilde\bold e_0=\cosh(\phi)\,\bold e_0+\sinh(\phi)\,\bold e_1,\\
\vspace{1ex}
&\tilde\bold e_1=\sinh(\phi)\,\bold e_0+\cosh(\phi)\,\bold e_1.
\endaligned
\mytag{2.2}
$$
Under the basis transformation \mythetag{2.2} the matrices 
\mythetag{2.1} are transformed according to the standard
tensorial rule:
$$
\xalignat 2
&\hskip -2em
\tilde g_{ij}=\sum^1_{p=0}\sum^1_{q=0}
S^p_i\,S^q_j\,g_{p\kern 0.4pt q},
&&\Check{\Tilde g}_{ij}=\sum^1_{p=0}\sum^1_{q=0}
S^p_i\,S^q_j\,\check g_{p\kern 0.4pt q},
\mytag{2.3}
\endxalignat
$$
where the components of the transition matrix $S$ are determined
by \mythetag{2.2}:
$$
S=\Vmatrix \cosh(\phi) &\sinh(\phi)\\ 
\vspace{2ex}
\sinh(\phi) & \cosh(\phi)\endVmatrix.
\mytag{2.4}
$$
Substituting \mythetag{2.4} into \mythetag{2.3}, one easily 
finds that the first matrix \mythetag{2.1} is invariant under
the Lorentz transformation \mythetag{2.2}, i\.\,e\.
$\tilde g_{ij}=g_{ij}$. For the non-diagonal matrix element 
of the second metric $\check\bold g$ in the new basis
$\tilde\bold e_0,\,\tilde\bold e_1$ we have
$$
\hskip -2em
\Check{\Tilde g}_{01}=\frac{\check g_{00}+\check g_{11}}{2}
\,\sinh(2\,\phi)+\check g_{01}\,\cosh(2\,\phi).
\mytag{2.5}
$$
The metric $\check\bold g$ is diagonalized simultaneously with
the metric $\bold g$ if the equation $\Check{\Tilde g}_{01}=0$
can be solved with respect to $\phi$. Due to the above formula 
\mythetag{2.5} the equation $\Check{\Tilde g}_{01}=0$ is equivalent 
to
$$
\hskip -2em
(\check g_{00}+\check g_{11})
\,\tanh(2\,\phi)=-2\,\check g_{01}.
\mytag{2.6}
$$
Looking at \mythetag{2.6}, we define the following four 
mutually exclusive cases:
$$
\vcenter{\hsize 10cm
\offinterlineskip\settabs\+\indent
\vtrule
\hskip 3.6cm &\vtrule 
\hskip 5.8cm &\vtrule 
\cr\hrule 
\+\vtrule
\hss\quad The first case\hss&\vtrule
\hss\quad $\check g_{00}+\check g_{11}=0$ \ and 
\ $\check g_{01}=0$\hss &\vtrule\cr\hrule
\+\vtrule
\hss\quad The second case\hss&\vtrule
\hss\quad $\check g_{00}+\check g_{11}\neq 0$ \ and 
\ $\check g_{01}=0$\hss&\vtrule\cr\hrule
\+\vtrule
\hss\quad The third case\hss&\vtrule
\hss\quad $\check g_{00}+\check g_{11}=0$ \ and 
\ $\check g_{01}\neq0$\hss&\vtrule\cr\hrule
\+\vtrule
\hss\quad The fourth case\hss&\vtrule
\hss\quad $\check g_{00}+\check g_{11}\neq 0$ \ and 
\ $\check g_{01}\neq 0$
\hss&\vtrule\cr\hrule}
\mytag{2.7}
$$\par
     In the {\bf first case} $\check g_{01}=0$. Therefore the second 
metric $\check\bold g$ is diagonal in the initial basis $\bold e_0,\,
\bold e_1$. Moreover, $\check g_{11}=-\check g_{00}$. If we denote
$\check g_{00}=a$, then \mythetag{2.1} is written as
$$
\xalignat 2
&\hskip -2em
g_{ij}=\Vmatrix\format\l\ &\ \r\\1 & 0\\ 0 & -1\endVmatrix,
&&\check g_{ij}=\Vmatrix\format\l\ &\ \r\\ a & 0\\
\vspace{1ex}0 & -a\endVmatrix.
\mytag{2.8}
\endxalignat
$$
As we see in \mythetag{2.8}, the metrics $\bold g$ and $\check\bold g$ 
in the first case do coincide up to the numeric factor $a$, i\.\,e\.
$\check\bold g=a\,\bold g$. They are always diagonalized simultaneously.
\par
     In the {\bf second case} $\check g_{01}=0$ too. Both metrics 
$\bold g$ and $\check\bold g$ are diagonal simultaneously in the initial 
basis $\bold e_0,\,\bold e_1$. If we denote $\check g_{00}=a$ and 
$\check g_{11}=b$, we get
$$
\xalignat 2
&\hskip -2em
g_{ij}=\Vmatrix\format\l\ &\ \r\\1 & 0\\ 0 & -1\endVmatrix,
&&\check g_{ij}=\Vmatrix\format\l\ &\ \r\\ a & 0\\
\vspace{1ex}0 & \ b\endVmatrix.
\mytag{2.9}
\endxalignat
$$
The equation \mythetag{2.6} in this case is solvable and its solution is
$\phi=0$. This means that \mythetag{2.2} is the identical transformation
where $\tilde\bold e_0=\bold e_0$ and $\tilde\bold e_1=\bold e_1$.\par
     In the {\bf third case} $\check g_{01}\neq 0$, i\.\,e\. the second
metric $\check\bold g$ is not diagonal, but we have the relationship 
$\check g_{11}=-\check g_{00}$. If we denote $\check g_{00}=a$ and 
$\check g_{01}=b$, then we get
$$
\xalignat 2
&\hskip -2em
g_{ij}=\Vmatrix\format\l\ &\ \r\\1 & 0\\ 0 & -1\endVmatrix,
&&\check g_{ij}=\Vmatrix\format\l\ &\ \r\\ a & b\\
\vspace{1ex}b & -a\endVmatrix.
\mytag{2.10}
\endxalignat
$$
The equation \mythetag{2.6} in this case is not solvable, i\.\,e\. the
metrics $\bold g$ and $\check\bold g$ cannot be diagonalized simultaneously.
Therefore we postulate \mythetag{2.10} to be another canonical presentation
for the metrics $\bold g$ and $\check\bold g$.\par
     The fourth case in the table \mythetag{2.7} is subdivided into three
subcases 
$$
\vcenter{\hsize 10cm
\offinterlineskip\settabs\+\indent
\vtrule
\hskip 3.6cm &\vtrule 
\hskip 5.8cm &\vtrule 
\cr\hrule 
\+\vtrule
\hss\quad Case 4, subcase 1\hss&\vtrule
\hss\quad $|2\,\check g_{01}|<|\check g_{00}+\check g_{11}|$
\hss &\vtrule\cr\hrule
\+\vtrule
\hss\quad Case 4, subcase 2\hss&\vtrule
\hss\quad $|2\,\check g_{01}|>|\check g_{00}+\check g_{11}|$
\hss&\vtrule\cr\hrule
\+\vtrule
\hss\quad Case 4, subcase 3\hss&\vtrule
\hss\quad $|2\,\check g_{01}|=|\check g_{00}+\check g_{11}|$
\hss&\vtrule\cr\hrule}
\mytag{2.11}
$$\par
     In the {\bf first subcase} of the case 4 the equation 
\mythetag{2.6} is solvable. Indeed, since $\check g_{00}
+\check g_{11}\neq 0$, we can write it as follows:
$$
\tanh(2\,\phi)=-\frac{2\,\check g_{01}}{\check g_{00}+\check g_{11}}.
\mytag{2.12}
$$
The function $\tanh(2\,\phi)$ is a growing smooth function on the
real axis $\Bbb R$, its values range from $-1$ as $\phi\to -\infty$
to $+1$ as $\phi\to +\infty$. For this reason the equation
\mythetag{2.12} has a unique solution $\phi=\phi_0$. Substituting
it into \mythetag{2.2}, we find a new basis $\tilde\bold e_0,\,
\tilde\bold e_1$ where both metrics $\bold g$ and $\check\bold g$ 
are diagonal. Their matrices take their canonical forms \mythetag{2.9} 
in this new basis.\par
     In the {\bf second subcase} of the case 4 the equation 
\mythetag{2.6} is not solvable. Therefore we take the sum
$\Check{\Tilde g}_{00}+\Check{\Tilde g}_{11}$. The vanishing condition
for this sum leads to the following equation for the parameter $\phi$:
$$
\hskip -2em
\tanh(2\,\phi)=-\frac{\check g_{00}+\check g_{11}}
{2\,\check g_{01}}.
\mytag{2.13}
$$
Looking at the second raw in the table \mythetag{2.11}, we see that
the second subcase of the case 4 is that very case where the equation
\mythetag{2.13} is solvable and has a unique solution $\phi=\phi_0$
Substituting this solution into \mythetag{2.2}, we find a new basis 
$\tilde\bold e_0,\,\tilde\bold e_1$ where the matrices of the metrics
$\bold g$ and $\check\bold g$ take their canonical forms \mythetag{2.10}.
\par
     The {\bf third subcase} of the case 4 is a special case. It subdivides 
into two subcases of the next level. They are listed in the following table:
$$
\vcenter{\hsize 10cm
\offinterlineskip\settabs\+\indent
\vtrule
\hskip 3.6cm &\vtrule 
\hskip 5.8cm &\vtrule 
\cr\hrule 
\+\vtrule
\hss\quad Case 4, subcase 3A\hss&\vtrule
\hss\quad $2\,\check g_{01}=\check g_{00}+\check g_{11}$
\hss &\vtrule\cr\hrule
\+\vtrule
\hss\quad Case 4, subcase 3B\hss&\vtrule
\hss\quad $-2\,\check g_{01}=\check g_{00}+\check g_{11}$
\hss&\vtrule\cr\hrule}
\mytag{2.14}
$$
These two subcases \mythetag{2.14} are studied in the next section.
\head
3. Associated operators.
\endhead
     The first metric $\bold g$ is a Minkowski type metric with the
signature $(+,-)$. It is non-degenerate. For this reason we can define
the associated operator $\check\bold F$ for the second metric with
respect to the first one. It is introduced by the formula
$$
\hskip -2em
g(\check\bold F(\bold X),\bold Y)=\check g(\bold X,\bold Y).
\mytag{3.1}
$$
Here $\bold X$ and $\bold Y$ are two arbitrary vectors of the space $V$.
Due to the symmetry of the quadratic forms $g$ and $\check g$ we
can extend \mythetag{3.1} as follows:
$$
\hskip -2em
g(\check\bold F(\bold X),\bold Y)=\check g(\bold X,\bold Y)
=g(\bold X,\check\bold F(\bold Y)),
\mytag{3.2}
$$
The formulas \mythetag{3.2} mean that $\check\bold  F$ is a symmetric 
operator with respect to the metric $\bold g$. In the coordinate form
the associated operator $\check\bold F$ is represented by a matrix:
$$
\hskip -2em
\check F^i_j=
\Vmatrix\check F^0_0 & \check F^0_1\\
\vspace{1ex}\check F^1_0 & \check F^1_1\endVmatrix.
\mytag{3.3}
$$
The components of the matrix \mythetag{3.3} are given by the formula
$$
\check F^i_j=\sum^3_{s=0}g^{is}\,\check g_{sj}.
$$
Here $g^{is}$ are the components of the matrix inverse to the matrix 
of the first metric $\bold g$. Applying this formula to \mythetag{2.1},
for $\check\bold F$ in the basis $\bold e_0,\,\bold e_1$ we get
$$
\hskip -2em
\check F^i_j=
\Vmatrix\check g_{00} & \check g_{01}\\
\vspace{1ex}-\check g_{01} & -\check g_{11}\endVmatrix.
\mytag{3.4}
$$
Using \mythetag{3.4}, we can calculate the invariants for the pair 
of metrics $\bold g$ and $\check\bold g$:
$$
\xalignat 2
&\hskip -2em
\tr\check\bold F=\check g_{00} -\check g_{11},
&&\det\check\bold F=(\check g_{01})^2-\check g_{00}\,\check g_{11}.
\mytag{3.5}
\endxalignat
$$
Relying on \mythetag{3.5}, we perform the following calculations:
$$
\hskip -2em
\gathered
(\check g_{00} +\check g_{11})^2-4\,(\check g_{01})^2
=(\check g_{00} -\check g_{11})^2\,+\\
\vspace{1ex}
+\,4\,\check g_{00}\,\check g_{11}
-4\,(\check g_{01})^2=(\tr\check\bold F)^2
-4\,\det\check\bold F.
\endgathered
\mytag{3.6}
$$
Due to the formula \mythetag{3.6} we can write the conditions in 
the table \mythetag{2.11} in the invariant coordinate-free form:
$$
\vcenter{\hsize 10cm
\offinterlineskip\settabs\+\indent
\vtrule
\hskip 3.6cm &\vtrule 
\hskip 5.8cm &\vtrule 
\cr\hrule 
\+\vtrule
\hss\quad Case 4, subcase 1\hss&\vtrule
\hss\quad $(\tr\check\bold F)^2>4\,\det\check\bold F$
\hss &\vtrule\cr\hrule
\+\vtrule
\hss\quad Case 4, subcase 2\hss&\vtrule
\hss\quad $(\tr\check\bold F)^2<4\,\det\check\bold F$
\hss&\vtrule\cr\hrule
\+\vtrule
\hss\quad Case 4, subcase 3\hss&\vtrule
\hss\quad $(\tr\check\bold F)^2=4\,\det\check\bold F$
\hss&\vtrule\cr\hrule}
\mytag{3.7}
$$\par
     In the subcase 1 of the table \mythetag{3.7} the associated 
operator $\check\bold F$ has two real eigenvalues $\lambda_0\neq
\lambda_1$. Let $\bold v_0$ and $\bold v_1$ be the eigenvectors
of the operator $\check\bold F$ corresponding to the eigenvalues
$\lambda_0$ and $\lambda_1$ respectively. Since $\lambda_0\neq
\lambda_1$, they are orthogonal to each other with respect to 
both metrics $\bold g$ and $\check\bold g$:
$$
\xalignat 2
&\hskip -2em
g(\bold v_0,\bold v_1)=0,
&&\check g(\bold v_0,\bold v_1)=0.
\mytag{3.8}
\endxalignat
$$
The proof of this fact is derived from \mythetag{3.2}. Indeed, we 
have
$$
\lambda_0\,g(\bold v_0,\bold v_1)
=g(\check\bold F(\bold v_0),\bold v_1)=\check g(\bold v_0,\bold v_1)
=g(\bold v_0,\check\bold F(\bold v_1))
=\lambda_1\,g(\bold v_0,\bold v_1).
\quad
\mytag{3.9}
$$
From \mythetag{3.9} we derive $(\lambda_1-\lambda_0)\,g(\bold v_0,
\bold v_1)=0$, which yields $g(\bold v_0,\bold v_1)=0$. Substituting
this equality back to the formulas \mythetag{3.9}, we get 
$\check g(\bold v_0,\bold v_1)$. Thus, both equalities \mythetag{3.8} 
are proved.\par
     If we choose the vectors $\bold v_0,\,\bold v_1$ for a basis, 
then the equalities \mythetag{3.9} mean that both metrics $\bold g$ 
and $\check\bold g$ are diagonal in this basis. The signature of the
metric $\bold g$ is $(+,-)$. For this reason $g(\bold v_0,\bold v_0)$
and $g(\bold v_1,\bold v_1)$ are two nonzero numbers of opposite signs.
Without loss of generality we can assume that $g(\bold v_0,\bold v_0)$
is positive and $g(\bold v_1,\bold v_1)$ is negative. We can 
normalize the eigenvectors $\bold v_0$ and $\bold v_1$ so that 
$g(\bold v_0,\bold v_0)=1$ and $g(\bold v_1,\bold v_1)=-1$. Then
$\bold v_0,\,\bold v_1$ is that very basis, where the metrics
$\bold g$ and $\check\bold g$ take their canonical forms 
\mythetag{2.9} with $a=\lambda_0$ and $b=-\lambda_1$.
\mytheorem{3.1} If\/ $(\tr\check\bold F)^2>4\,\det\check\bold F$, then
the metrics $\bold g$ and $\check\bold g$ are given by the matrices
\mythetag{2.9} with $a+b\neq 0$ in a basis composed by eigenvectors 
of the operator $\check\bold F$.
\endproclaim 
     In the subcase 2 of the table \mythetag{3.7} the associated 
operator $\check\bold F$ has two complex eigenvalues conjugate to
each other: $\lambda_1=\overline{\,\lambda_0}$. More exactly, 
$\lambda_0\neq\lambda_1$ are the eigenvalues of the complexified
operator $\check\bold F$ in the complexification $\Bbb CV=\Bbb C\otimes
V$ of the vector space $V$. The complex space $\Bbb CV$ is naturally
equipped with the involution of complex conjugation:
$$
\hskip -2em
\tau\!:\,\Bbb CV\to\Bbb CV.
\mytag{3.10}
$$
The space $V$ is embedded into $\Bbb CV$ as a $\Bbb R$-linear subspace 
invariant under the involution \mythetag{3.10}. Since $\check\bold F$
is a complexification of an operator acting in $V$, it commutes with
$\tau$. Therefore, if $\bold v_0$ is an eigenvector corresponding to
the eigenvalue $\lambda_0$, then $\bold v_1=\tau(\bold v_0)$ is an
eigenvector corresponding to the eigenvalue 
$\lambda_1=\overline{\,\lambda_0}$. Let's define the following two 
vectors:
$$
\xalignat 2
&\bold e_0=\frac{\bold v_0+\bold v_1}{\sqrt{2}}
&&\bold e_1=\frac{\bold v_0-\bold v_1}{\sqrt{2}\,i}.
\mytag{3.11}
\endxalignat
$$
The vectors \mythetag{3.11} are invariant under the action of the
involution $\tau$. Hence they belong to $V$. These vectors are nonzero
and linearly independent. They form a basis in $V$. Applying $\check
\bold F$ to \mythetag{3.11}, we find
$$
\hskip -2em
\aligned
&\check\bold F(\bold e_0)=\frac{\lambda_0\,\bold v_0
+\overline{\,\lambda_0}\,\bold v_1}{\sqrt{2}}=\Rea(\lambda_0)
\,\bold e_0-\Img(\lambda_0)\,\bold e_1,\\
&\check\bold F(\bold e_1)=\frac{\lambda_0\,\bold v_0
-\overline{\,\lambda_0}\,\bold v_1}{\sqrt{2}\,i}=\Img(\lambda_0)
\,\bold e_0+\Rea(\lambda_0)\,\bold e_1.
\endaligned
\mytag{3.12}
$$
Note that the vectors $\bold v_0$ and $\bold v_1$ are orthogonal to each 
other with respect to both metrics $\bold g$ and $\check\bold g$, i\.\,e\.
the formulas \mythetag{3.8} are valid. The arguments for that here are the 
same as in \mythetag{3.9}. Due to \mythetag{3.8} the metric $\bold g$
is diagonal in the basis $\bold v_0,\,\bold v_1$. It is a non-degenerate
metric. Hence, $g(\bold v_0,\bold v_0)$ and $g(\bold v_1,\bold v_1)$ are
nonzero. Due to the complexity of the space $\Bbb CV$ the vectors 
$\bold v_0$ and $\bold v_1$ can be normalized to the unity:
$$
\pagebreak
\xalignat 2
&\hskip -2em
g(\bold v_0,\bold v_0)=1,
&&g(\bold v_1,\bold v_1)=1.
\mytag{3.13}
\endxalignat
$$
From \mythetag{3.11}, \mythetag{3.8} and \mythetag{3.13} we easily derive
$$
\xalignat 3
&\hskip -2em
g(\bold e_0,\bold e_0)=1,
&&g(\bold e_0,\bold e_1)=0,
&&g(\bold e_1,\bold e_1)=-1.
\mytag{3.14}
\endxalignat
$$
Let's denote $\Rea(\lambda_0)=a$ and $\Img(\lambda_0)=b$. Then from 
\mythetag{3.12} and \mythetag{3.14}, using the formula \mythetag{3.2},
we derive that the metrics $\bold g$ and $\check\bold g$ take their
canonical forms \mythetag{2.10}.
\mytheorem{3.2} If\/ $(\tr\check\bold F)^2<4\,\det\check\bold F$, 
then the metrics $\bold g$ and $\check\bold g$ are given by the 
matrices \mythetag{2.10} in a basis produced from eigenvectors of 
the complexified associated operator $\check\bold F$ according to
the formulas \mythetag{3.11}.
\endproclaim 
     Now let's proceed to the subcase 3 of the table \mythetag{3.7}. 
In this case the associated operator $\check\bold F$ has one real 
eigenvalue $\lambda_0$ of the multiplicity $2$. Assume that the
metrics $\bold g$ and $\check\bold g$ are brought to the form 
\mythetag{2.1} in some basis $\bold e_0,\bold e_1$. Then this
subcase 3 is subdivided into two subcases 3A and 3B of the next
level (see the table \mythetag{2.14}). Actually, the subcase 3B 
is equivalent to the subcase 3A. Indeed, assume that the
condition $-2\,\check g_{01}=\check g_{00}+\check g_{11}$ is
fulfilled in the basis $\bold e_0,\bold e_1$. Then we perform
the following basis transformation:
$$
\xalignat 2
&\hskip -2em
\tilde\bold e_0=\bold e_0,
&&\tilde\bold e_1=-\bold e_1.
\mytag{3.15}
\endxalignat
$$
The transformation \mythetag{3.15} is characterized by the diagonal 
transition matrix
$$
\hskip -2em
S=\Vmatrix\format\l\ &\ \r\\1 & 0\\ 0 & -1\endVmatrix.
\mytag{3.16}
$$
Substituting \mythetag{3.16} into \mythetag{2.3}, we find that
$\tilde g_{ij}=g_{ij}$, i\.\,e\. the matrix of the metric $\bold g$ 
is invariant under the basis transformation \mythetag{3.15}, while 
for the matrices of the second metric we have the following
relationships:
$$
\xalignat 3
&\Check{\Tilde g}_{00}=g_{00},
&&\Check{\Tilde g}_{01}=-g_{01},
&&\Check{\Tilde g}_{11}=g_{11}.
\quad
\mytag{3.17}
\endxalignat
$$
Due to the formulas \mythetag{3.17}, from $-2\,\check g_{01}
=\check g_{00}+\check g_{11}$ we derive $2\,\Check{\Tilde g}_{01}
=\Check{\Tilde g}_{00}+\Check{\Tilde g}_{11}$. Thus, the subcase 
3B occurring in some basis $\bold e_0,\bold e_1$ can be transformed 
to the subcase 3A in some other basis.\par
     Continuing the study of the subcase 3 in \mythetag{3.7}, we 
restrict ourselves to the subcase 3A. Using the equality 
$2\,\check g_{01}=\check g_{00}+\check g_{11}$ we express 
$\check g_{01}$ through $\check g_{00}$ and $\check g_{11}$:
$$
\hskip -2em
\check g_{01}=\frac{\check g_{00}+\check g_{11}}{2}.
\mytag{3.18}
$$
Then we substitute the expression \mythetag{3.18} into the matrix 
\mythetag{3.4} and calculate the eigenvalue of the associated 
operator $\check\bold F$:
$$
\lambda_0=\frac{\check g_{00}-\check g_{11}}{2}.
\mytag{3.19}
$$
As we mentioned above, the operator $\check\bold F$ in this case 
has exactly one eigenvalue \mythetag{3.19} of the multiplicity $2$.
Let $A=\check F-\lambda_0\,I$, where $\check F$ is the matrix of the
operator $\check\bold F$ and $I$ is the unit matrix. Then, using
\mythetag{3.19} for $\lambda_0$, we obtain
$$
\hskip -2em
A^i_j=\frac{\check g_{00}+\check g_{11}}{2}\,
\Vmatrix\format\r\ &\ \r\\ 1 & 1\\
\vspace{1ex} -1 & -1\endVmatrix.
\mytag{3.20}
$$
Now, relying on \mythetag{3.20} we define the following quantity:
$$
\hskip -2em
\sigma=\sign(\check g_{00}+\check g_{11})=\cases
+1&\text{if \ \ }\check g_{00}+\check g_{11}>0;\\
\ \ 0&\text{if \ \ }\check g_{00}+\check g_{11}=0;\\
-1&\text{if \ \ }\check g_{00}+\check g_{11}<0.
\endcases
\mytag{3.21}
$$
We subdivide the subcase 3 in \mythetag{3.7} into three subcases 
of the next level regarding the value of $\sigma$ in \mythetag{3.21}. 
They are listed in the table
$$
\vcenter{\hsize 10cm
\offinterlineskip\settabs\+\indent
\vtrule
\hskip 3.6cm &\vtrule 
\hskip 3cm &\vtrule 
\cr\hrule 
\+\vtrule
\hss\quad Case 4, subcase 3(1)\hss&\vtrule
\hss\quad $\sigma=1$
\hss &\vtrule\cr\hrule
\+\vtrule
\hss\quad Case 4, subcase 3(2)\hss&\vtrule
\hss\quad $\sigma=-1$
\hss&\vtrule\cr\hrule
\+\vtrule
\hss\quad Case 4, subcase 3(3)\hss&\vtrule
\hss\quad $\sigma=0$
\hss&\vtrule\cr\hrule}
\mytag{3.22}
$$\par
     The subcase 3(3) is the most simple in the table \mythetag{3.22}.
In this case the matrix \mythetag{3.20} is equal to zero, i\.\,e\.
$g_{00}+g_{11}=0$. Then we denote
$$
\hskip -2em
\check g_{00}=-\check g_{11}=a.
\mytag{3.23}
$$
Substituting \mythetag{3.23} into \mythetag{3.18} and \mythetag{3.19}, 
we find that
$$
\xalignat 2
&\hskip -2em
\lambda_0=a, &&\check g_{01}=0.
\mytag{3.24}
\endxalignat
$$
Substituting \mythetag{3.23} and \mythetag{3.24} back into \mythetag{2.1},
we see that the subcase 3(3) in the table \mythetag{3.22} is equivalent
to the first subcase in the table \mythetag{2.7}. 
\mytheorem{3.3}If\/ $(\tr\check\bold F)^2=4\,\det\check\bold F$ and 
$\sigma=0$, then the metrics $\bold g$ and $\check\bold g$ differ only 
by a scalar factor. They can be brought to the canonical form 
\mythetag{2.8} in some \nolinebreak basis.
\endproclaim
     Let's proceed to the subcase 3(1) in the table \mythetag{3.22}. In 
this case $\check g_{00}+\check g_{11}>0$. Therefore we denote $\check g_{00}
+\check g_{11}=2\,\beta^2$ and $\lambda_0=a$. Then \mythetag{3.18} and 
\mythetag{3.19} yield
$$
\xalignat 3
&\hskip -2em
\check g_{00}=\beta^2+a, 
&&\check g_{01}=\beta^2,
&&\check g_{11}=\beta^2-a.
\qquad
\mytag{3.25}
\endxalignat
$$
Substituting \mythetag{3.25} into \mythetag{2.1}, we find
$$
\xalignat 2
&\hskip -2em
g_{ij}=\Vmatrix\format\l\ &\ \r\\1 & 0\\ 
\vspace{1ex}0 & -1\endVmatrix,
&&\check g_{ij}=\Vmatrix\beta^2+a & \ \ \beta^2\\
\vspace{1ex}\vspace{1ex}\beta^2 & \beta^2-a\endVmatrix.
\mytag{3.26}
\endxalignat
$$
The matrices \mythetag{3.26} present the metrics $\bold g$ and 
$\check\bold g$ in some basis $\bold e_0,\,\bold e_1$. Now we 
perform the following basis transformation:
$$
\xalignat 2
&\hskip -2em
\tilde\bold e_1=\frac{1}{2\,\beta}\,\bold e_0
+\frac{1}{2\,\beta}\,\bold e_1,
&&\tilde\bold e_0=\beta\,\bold e_0-\beta\,\bold e_1.
\mytag{3.27}
\endxalignat
$$
Upon performing the basis transformation \mythetag{3.27} we find that 
the metrics $\bold g$ and $\check\bold g$ are presented by the matrices
$$
\pagebreak
\xalignat 2
&\hskip -2em
g_{ij}=\Vmatrix\format\l\ &\ \r\\0 & 1\\ 
1 & 0\endVmatrix,
&&\check g_{ij}=\Vmatrix\format\l\ &\ \r\\
1 & a\\a & 0\endVmatrix
\mytag{3.28}
\endxalignat
$$
in the new basis. The presentation \mythetag{3.28} is a canonical
presentation for the metric pair $\bold g$, $\check\bold g$ in the
subcase 3(1).
\mytheorem{3.4}If\/ $(\tr\check\bold F)^2=4\,\det\check\bold F$ and 
$\sigma=1$, then the metrics $\bold g$ and $\check\bold g$ are presented
by the matrices \mythetag{3.28} in some basis.
\endproclaim
     The subcase 3(2) is similar to the subcase 3(1). In this case
$\check g_{00}+\check g_{11}<0$ Therefore we denote $\check g_{00}
+\check g_{11}=-2\,\beta^2$ and $\lambda_0=a$. Then \mythetag{3.18} 
and \mythetag{3.19} yield
$$
\xalignat 3
&\hskip -2em
\check g_{00}=a-\beta^2, 
&&\check g_{01}=-\beta^2,
&&\check g_{11}=-a-\beta^2.
\qquad
\mytag{3.29}
\endxalignat
$$
Due to \mythetag{3.29} the formulas \mythetag{2.1} specialize to the 
following ones:
$$
\xalignat 2
&\hskip -2em
g_{ij}=\Vmatrix\format\l\ &\ \r\\1 & 0\\ 
\vspace{1ex}0 & -1\endVmatrix,
&&\check g_{ij}=\Vmatrix a-\beta^2 & \ \ -\beta^2\\
\vspace{1ex}\vspace{1ex}-\beta^2 & -a-\beta^2\endVmatrix.
\mytag{3.30}
\endxalignat
$$
Now we perform the following basis transformation:
$$
\xalignat 2
&\hskip -2em
\tilde\bold e_0=\beta\,\bold e_0-\beta\,\bold e_1,
&&\tilde\bold e_1=\frac{1}{2\,\beta}\,\bold e_0
+\frac{1}{2\,\beta}\,\bold e_1.
\mytag{3.31}
\endxalignat
$$
By means of \mythetag{3.31} we bring the matrices \mythetag{3.30} 
to their canonical forms:
$$
\xalignat 2
&\hskip -2em
g_{ij}=\Vmatrix\format\l\ &\ \r\\0 & 1\\ 
1 & 0\endVmatrix,
&&\check g_{ij}=\Vmatrix\format\l\ &\ \r\\
 0 & a\\a & -1\endVmatrix.
\mytag{3.32}
\endxalignat
$$
\mytheorem{3.5}If\/ $(\tr\check\bold F)^2=4\,\det\check\bold F$ and 
$\sigma=-1$, then the metrics $\bold g$ and $\check\bold g$ are presented
by the matrices \mythetag{3.32} in some basis.
\endproclaim
\head
4. Classification.
\endhead
     The cases and subcases considered in the previous two sections 
are excessive. Some of them are equivalent to others and some of them 
are particular cases of others. The actual classification of metric 
pairs, one of which is a Minkowski type metric, is given by the 
theorems~\mythetheorem{3.1}, \mythetheorem{3.2}, \mythetheorem{3.3}, 
\mythetheorem{3.4}, and \mythetheorem{3.5}. We gather the results of 
these theorems into the following table:
$$
\vcenter{\hsize=200pt
\offinterlineskip
\halign{\vrule#&\hfill\quad #\quad\hfill&\vrule#&\quad #\quad\hfill
&\vrule#\cr
\noalign{\hrule}
height 14pt depth 8pt & Condition &&\hfill Canonical presentation&\cr
\noalign{\hrule}
height 18pt depth 12pt &$(\tr\check\bold F)^2>4\,\det\check\bold F$ 
&&$g_{ij}=\Vmatrix\format\l\ &\ \r\\1 & 0\\ 0 & -1\endVmatrix$,\quad
$\check g_{ij}=\Vmatrix\format\l\ &\ \r\\ a & 0\\
\vspace{1ex}0 & \ b\endVmatrix$ with $b\neq -a$&\cr
\noalign{\hrule}
height 18pt depth 12pt &$(\tr\check\bold F)^2<4\,\det\check\bold F$ 
&&$g_{ij}=\Vmatrix\format\l\ &\ \r\\1 & 0\\ 0 & -1\endVmatrix$,\quad
$\check g_{ij}=\Vmatrix\format\l\ &\ \r\\ a & b\\
\vspace{1ex}b & -a\endVmatrix$ with $b\neq 0$&\cr
\noalign{\hrule}
height 18pt depth 12pt &$\vcenter{\hsize=35em\noindent 
$(\tr\check\bold F)^2=4\,\det\check\bold F$\newline\vphantom{!}\quad 
\ and \strut $\sigma=0$}$\kern -27.5em
&&$g_{ij}=\Vmatrix\format\l\ &\ \r\\1 & 0\\ 0 & -1\endVmatrix$,\quad
$\check g_{ij}=\Vmatrix\format\l\ &\ \r\\ a & 0\\
\vspace{1ex}0 & -a\endVmatrix$&\cr
\noalign{\hrule}
height 18pt depth 12pt &$\vcenter{\hsize=35em\noindent 
$(\tr\check\bold F)^2=4\,\det\check\bold F$\newline\vphantom{!}\quad 
\ and \strut $\sigma=1$}$\kern -27.5em
&&$g_{ij}=\Vmatrix\format\l\ &\ \r\\0 & 1\\ 1 & 0\endVmatrix$,\quad
\ \ $\check g_{ij}=\Vmatrix\format\l\ &\ \r\\ 1 & a\\a & 0\endVmatrix$&\cr
\noalign{\hrule}
height 18pt depth 12pt &$\vcenter{\hsize=35em\noindent 
$(\tr\check\bold F)^2=4\,\det\check\bold F$\newline\vphantom{!}\quad 
\ and \strut $\sigma=-1$}$\kern -27.5em
&&$g_{ij}=\Vmatrix\format\l\ &\ \r\\0 & 1\\ 1 & 0\endVmatrix$,\quad
\ \ $\check g_{ij}=\Vmatrix\format\l\ &\ \r\\ 0 & a\\a & -1\endVmatrix$
&\cr
\noalign{\hrule}}}
\mytag{4.1}
$$
The quantity $\sigma$ in the table \mythetag{4.1} is an invariant
of a pair of metrics. It is very important to note that this
invariant cannot be expressed through the invariants of the
associated operator (\,$\tr\check\bold F$ \,and \ $\det\check
\bold F$\,). The formula \mythetag{3.21} defines this invariant
in a special basis, where the first metric $\bold g$ is 
diagonalized:
$$
g_{ij}=\Vmatrix\format\l\ &\ \r\\1 & 0\\ 0 & -1\endVmatrix.
$$
However, there must be a formula or an algorithm for calculating 
the invariant $\sigma$ in an arbitrary basis without diagonalizing
the metric $\bold g$. 
\head
5. Dedicatory.
\endhead
     This paper is dedicated to my uncle Amir Minivalievich Nagaev.

\enddocument
\end